\documentclass[11pt]{article}
\usepackage{amssymb,amsmath,graphicx,theorem}
\usepackage{hyperref}

\def\GG{{\cal G}}

\def\AA{{\cal A}}

\def\R{{\mathbf R}}

\def\T{{^{\sf T}}}

\def\Aut{{\rm Aut}}
\def\orb{{\rm orb}}
\def\hom{\hbox{\rm hom}}

\def\rk{\hbox{\rm rk}}
\def\tr{\hbox{\rm tr}}

\newtheorem{theorem}{Theorem}[section]

\newtheorem{lemma}[theorem]{Lemma}
\newtheorem{claim}{Claim}[section]
\newtheorem{corollary}[theorem]{Corollary}
\theorembodyfont{\rmfamily}

\long\def\killtext#1{}

\newenvironment{proof}{\noindent{\bf Proof. }}{\hfill$\square$\medskip}

\begin{document}

\title{The rank of connection matrices and the dimension of graph algebras}
\author{{\sc L\'aszl\'o Lov\'asz}\\
Microsoft Research \\
One Microsoft Way\\
Redmond, WA 98052\\}

\date{August 2004\\
Microsoft Research Technical Report TR-2004-82}

\maketitle

\tableofcontents

\begin{abstract}
Connection matrices were introduced in \cite{FLS}, where they were
used to characterize graph homomorphism functions. The goal of this
note is to determine the exact rank of these matrices. The result can
be rephrased in terms of graph algebras, also introduced in
\cite{FLS}. Yet another version proves that if two $k$-tuples of
nodes behave the same way from the point of view of graph
homomorphisms, then they are equivalent under the automorphism group.
\end{abstract}

\section{Introduction}

For two finite graphs $F$ and $G$, let $\hom(F,G)$ denote the number
of homomorphisms (adjacency-preserving mappings) from $F$ to $G$.

For every fixed $G$, let us construct the following (infinite) matrix
$M(k,G)$. The rows and columns are indexed by finite graphs $F$ in
which $k$ nodes are labeled $1,\dots,k$ (there can be any number of
unlabeled nodes). The entry in the intersection of the row
corresponding to $F_1$ and the column corresponding to $F_2$ is
$\hom(F_1F_2,G)$, where $F_1F_2$ is obtained by considering the
disjoint union of $F_1$ and $F_2$, and identifying the nodes labeled
the same way. This matrix is called the {\it $k$-th connection matrix
for homomorphisms into $G$}.

One can extend this definition to the case when $G$ has edgeweights
and nodeweights (see Section \ref{HOMDEFS} for the exact
definitions). Connection matrices were introduced by Freedman,
Lov\'asz and Schrijver \cite{FLS}, where they were used to
characterize graph homomorphism functions $\hom(.,G)$. In particular,
it was shown that {\it connection matrices are positive semidefinite
and $M(k,G)$ has rank at most $|V(G)|^k$.} (We'll reproduce the
simple proof of this assertion in section \ref{HOMDEFS}. The main
result in \cite{FLS} is a converse to this statement, which we don't
quote here.)

This assertion raises the question: what is the exact rank of
$M(k,G)$? The aim of this paper is to determine this rank.

The operation of gluing together two graphs along their labeled nodes
gives rise to a commutative algebra defined on formal linear
combinations of graphs. This is a tool that was introduced in
\cite{FLS}, and will be very useful for us too. The results of this
paper can also be viewed as describing the dimension of these
algebras.

A third version of these results is motivated by the following. One
often classifies nodes of a graph by their degrees. We can also
consider the following stronger classification: for every simple
graph $F$ with a specified node, consider the number $\hom_v(F,G)$ of
those homomorphisms of $F$ into $G$ that map the specified node onto
$v$. This way each node $v\in V(G)$ is assigned an infinite vector
$h_v=(\hom_v(F_1,G),\hom_v(F_2,G),\dots)$, where $(F_1,F_2,\dots)$ is
any enumeration of all simple graphs with a specified node. Are there
any linear relations between these vectors? Clearly two vectors
$h_u,h_v$ are the same if there is an automorphism of $G$ that moves
$u$ to $v$. For unweighted graphs, this turns out to be all; for
weighted graphs, the situation is a bit more complicated, but we'll
determine all relations; they are still trivial in some sense. These
results extend to graphs $F$ with $k$ specified nodes instead of $1$.

These results have various applications; for example, Lov\'asz and
S\'os \cite{LS} use it to characterize generalized quasirandom
graphs.

\section{Homomorphisms and connection matrices}\label{HOMDEFS}

We start with extending the notions introduced above to weighted
graphs. A {\it weighted graph} $G$ is a graph with a positive real
weight $\alpha_G(i)$ associated with each node and a real weight
$\beta_G(i,j)$ associated with each edge $ij$. An edge with weight 0
will play the same role as no edge between those nodes, so we can
assume that all the edge weights are nonzero, or that $G$ is a
complete graph with loops at each node, whichever is more convenient.

Let $F$ be an unweighted graph (possibly with multiple edges, but no
loops) and $G$, a weighted graph. To every $\phi:~V(F)\to V(G)$, we
assign two weights:
\[
\alpha_\phi=\prod_{u\in V(F)} \alpha_G(\phi(u))
\]
and
\[
\hom_\phi(F,G)=\prod_{u,v\in V(F)} \beta_G(\phi(u),\phi(v)).
\]
Define
\[
\hom(F,G)=\sum_{\phi:~V(F)\to V(G)} \alpha_\phi\hom_\phi(F,G).
\]
If all the node-weights and edge-weights in $G$ are $1$, then this is
the number of homomorphisms from $F$ into $G$ (with no weights).

For the purpose of this paper, it will be convenient to assume that
$G$ is a complete graph with a loop at all nodes (missing edges can
be added with weight $0$). Then the weighted graph $G$ is completely
described by a positive real vector
$a=(\alpha_1,\dots,\alpha_m)\in\R^m$ and a real symmetric matrix
$B=(\beta_{ij})\in\R^{m\times m}$. It will be convenient to assume
that
\[
\sum_{i=1}^m \alpha_i=1;
\]
this only means scaling of the $\hom$ function by an appropriate
power of $\sum_i\alpha_i$, and will not influence the results.

A {\it $k$-labeled graph} ($k\ge 0$) is a finite graph in which $k$
nodes are labeled by $1,2,\dots k$. Two $k$-labeled graphs are {\it
isomorphic}, if there is a label-preserving isomorphism between them.
We denote by $K_k$ the $k$-labeled complete graph on $k$-nodes, and
by $E_k$, the $k$-labeled graph on $k$ nodes with no edges.

Let $F_1$ and $F_2$ be two $k$-labeled graphs. Their {\it product}
$F_1F_2$ is defined as follows: we take their disjoint union, and
then identify nodes with the same label. For two 0-labeled graphs,
$F_1F_2$ is just their disjoint union.

The definition of connection matrices can be extended to the case
when $G$ is weighted in a trivial way: The rows and columns of
$M(k,G)$ are indexed by isomorphism types of $k$-labeled graphs. The
entry in the intersection of the row corresponding to $F_1$ and the
column corresponding to $F_2$ is $\hom(F_1F_2,G)$. Let us also recall
their main properties:

\begin{lemma}\label{FLS-EASY}
The connection matrices $M(k,G)$ are positive semidefinite and
$M(k,G)$ has rank at most $|V(G)|^k$.
\end{lemma}

This lemma will follow very easily if we introduce two further
matrices. Let us extend our notation by defining, for any $k$-labeled
graph $F$ and mapping $\phi:\,[1,k]\to V(G)$,
\begin{equation}\label{PHI-HOM}
\hom_\phi(F,G)=\sum_{\psi:~V(F)\to V(G)\atop \psi
\text{~extends~}\phi}\frac{\alpha_\psi}{\alpha_\phi}\hom_\psi(F,G)
\end{equation}
So
\[
\hom(F,G)=\sum_{\phi:\,[1,k]\to V(G)} \alpha_\phi\hom_\phi(F,G).
\]
Furthermore, for any two $k$-labeled graph $F_1$ and $F_2$, we have
the important equation
\begin{equation}\label{PHI-PROD}
\hom_\phi(F_1F_2,G)=\hom_\phi(F_1,G)\hom_\phi(F_2,G).
\end{equation}
This expresses that once we mapped the common part, the mapping can
be extended to the rest of $F_1$ and $F_2$ independently.

Let $N(k,G)$ denote the matrix in which rows are indexed by maps
$\phi:~[1,k]\to V(G)$, columns are indexed by $k$-labeled graphs $F$,
and the entry in the intersection of the row $\phi$ and column $F$ is
$\hom_\phi(F,G)$. Let $A(k,G)$ denote the diagonal matrix whose rows
and columns are indexed by maps $\phi:~[1,k]\to V(G)$, and the
diagonal entry in row $\phi$ is $\alpha_\phi$. Then
(\ref{PHI-HOM})and (\ref{PHI-PROD}) imply
\begin{equation}\label{MN}
M(k,G)=N(k,G)\T A(k,G) N(k,G).
\end{equation}
This equation immediately implies that $M(k,G)$ is positive
semidefinite and
\[
\rk(M(k,G))=\rk(N(k,G))\le |V(G)|^k.
\]

When does the rank of a connection matrix attain this upper bound?
There are two (related, but different) types of degeneracy that
causes lower rank.

\medskip

\noindent{\bf Twins.} The first of these causes is easy to handle. We
call two nodes $i,j\in V(G)$ {\it twins}, if for every node $l\in
V(G)$, $\beta_{il}=\beta_{jl}$ (note: the condition includes $l=i$
and $l=j$; the node weights $\alpha_i$ play no role in this
definition). We say that $G$ is {\it twin-free}, if no two different
nodes are twins.

Suppose that $G$ is not twin-free, so that it has two twin nodes $i$
and $j$. Let us identify the equivalence classes of twin nodes,
define the node-weight $\alpha$ of a new node as the sum of the
node-weights of its pre-images, and define the weight of an edge as
the weight of any of its pre-images (which all have the same weight).
This way we get a twin-free graph $\bar{G}$ such that
$\hom(F,G)=\hom(F,\bar{G})$ for every graph $F$. It follows that the
rank of the connection matrices $M(k,G)$ and $M(k,\bar{G})$ are the
same, and this rank is at most $|V(\bar{G})|^k<|V(G)|^k$.

From now on, we assume that $G$ is twin-free.

\medskip

\noindent{\bf Automorphisms.} The second reason for rank loss in the
connection matrices will take more work to handle. Suppose that $G$
has a proper automorphism (a permutation of the nodes that preserves
both the node-weights and edge-weights). Then any two rows of
$N(k,G)$ defined by a mappings $\phi:~[1,k]\to V(G)$ and $\phi\sigma$
($\sigma\in\Aut(G)$) are equal. So the rank of $N(k,G)$ (and
$M(k,G)$) is at most the number of orbits of the automorphism group
of $G$ on ordered $k$-tuples of its nodes. The main result of this
paper is that equality holds here.

\begin{theorem}\label{AUTOMORPH}
Let $G$ be twin-free weighted graph. Let $\orb_k(G)$ denote the
number of orbits of the automorphism group of $G$ on ordered
$k$-tuples of its nodes. Then $\rk(M(k,G))=\orb_k(G)$ for every $k$.
\end{theorem}

\begin{corollary}\label{HOWGOOD}
Let $G$ be a weighted graph that has no twins and no automorphisms.
Then $\rk(M(k,G))=|V(G)|^k$ for every $k$.
\end{corollary}

Note that swapping twins $i$ and $j$ is almost an automorphism: the
only additional condition needed is that $\alpha_i=\alpha_j$. So in
particular, for unweighted graphs the condition that there are no
automorphisms implies that there are no twins.

Along the lines, we'll prove two lemmas, which are of independent
interest:

\begin{lemma}\label{HOMEQ}
Let $G$ be a twin-free weighted graph, let $\phi,\psi\in V(G)^k$, and
suppose that for every $k$-labeled graph $F$,
$\hom_\phi(F,G)=\hom_\psi(F,G)$. Then there exists an automorphism
$\sigma$ of $H$ such that $\psi=\phi\sigma$.
\end{lemma}

Fix an integer $k\ge 1$ and a weighted graph $G$. We say that a
vector $f:~V(G)^k\to\R$ is {\it invariant under automorphisms of
$G$}, if $f(\phi\sigma)=f(\phi)$ for every $\sigma\in\Aut(G)$.
Trivially, every column of $N(h,G)$ is invariant under automorphisms.

\begin{lemma}\label{HOMREP}
The column space of $N(k,G)$ consists of precisely those vectors
$f:~V(G)^k\to\R$ that are invariant under automorphisms of $G$.
\end{lemma}

As a final application, we prove an extension of an old result from
\cite{Lo} to weighted graphs:

\begin{corollary}\label{HOMDET}
Let $G_1$ and $G_2$ be twin-free weighted graphs, and assume that for
every simple graph $F$, $\hom(F,G_1)=\hom(F,G_2)$. Then $G_1$ and
$G_2$ are isomorphic.
\end{corollary}

\section{The algebra of graphs}\label{GRAPHALG}

A {\it $k$-labeled quantum graph} is a formal linear combination
(with real coefficients) of $k$-labeled graphs. Let $\GG_k$ denote
the (infinite dimensional) vector space of all $k$-labeled quantum
graphs. We can turn $\GG_k$ into an algebra by using $F_1F_2$
introduced above as the product of two generators, and then extending
this multiplication to the other elements linearly. Clearly $\GG_k$
is associative and commutative, and the empty graph $E_k$ is a unit
element in $\GG_k$.

We need to introduce some further (rather trivial) algebras. Let
$\AA_k$ be the algebra of formal linear combinations of maps
$\phi:\,[1,k]\to V(G)$, where multiplication is defined in a trivial
way: for two maps $\phi$ and $\psi$, let $\phi*\psi=\phi$ if
$\phi=\psi$ and 0 otherwise. The sum
\[
u_k=\sum_{\phi:\,[1,k]\to V(G)} \phi
\]
is the unit element of this algebra.

Next define $f(.)=\hom(.,G)$. Extend $f$ linearly to quantum graphs.
This function $f$ gives rise to additional structure. We introduce an
inner product on $\GG$ by
\begin{equation}\label{INNERPROD}
\langle x,y\rangle = f(xy).
\end{equation}
We'll see that this inner product is semidefinite: $\langle
x,x\rangle\ge 0$ for all $x$. We also introduce an inner product on
$\AA_k$ as follows: for two basis elements $\phi,\psi:~[1,k]\to
V(G)$, let
\[
\langle \phi,\psi\rangle =
  \begin{cases}
    \alpha_\phi & \text{if $\phi=\psi$}, \\
    0 & \text{otherwise},
  \end{cases}
\]
and then extend this bilinearly. Trivially, this inner product is
positive definite.

The function $f$ is multiplicative over connected components: if
$F_1,F_2\in\GG_0$, then
\[
f(F_1F_2)=f(F_1)f(F_2).
\]
This means that as a map $\GG_0\to\AA_{0}$ it an algebra
homomorphism.

The graph $G$ also gives rise to a map $f_k:~\GG_k\to\AA_k$ by
\[
f_k(F)=\sum_{\phi:~[1,k]\to[1,m]} \hom_\phi(F,G) \phi.
\]
We extend this map linearly to quantum graphs. For two $k$-labeled
quantum graphs $x,y$ we say that
\[
x\equiv y \pmod G
\]
if $f_k(x)=f_k(y)$.

It is easy to check that the mapping $f_k:~\GG_k\to\AA_k$ is an
algebra homomorphism, and preserves inner product. This in particular
implies that the inner product defined by (\ref{INNERPROD}) is
positive semidefinite. Furthermore, since the inner product in
$\AA_{k}$ is positive definite, the kernel of $f_k$ is exactly the
nullspace of the inner product (\ref{INNERPROD}). If we factor out
this nullspace, we get an algebra $\GG_k'$. It is easy to check that
\begin{equation}\label{DIMRANK}
\dim(\GG'_k)=\rk(M(k,G)).
\end{equation}
So we know that this dimension is at most $|V(G)|^k$; in particular
it is finite.

For every $k>0$, we define the {\it trace} $\tr:~\GG_k\to\GG_{k-1}$
simply by erasing the label $k$. We also have a linear map
$\tr:~\AA_k\to\AA_{k-1}$ defined by
\[
\tr(e_{i_1}\otimes \dots\otimes e_{i_k})=\alpha_{i_k}(e_{i_1}\otimes
\dots\otimes e_{i_{k-1}}).
\]
These operators correspond to each other in the sense that for every
$k$-labeled graph $F$,
\[
\tr(f_k(F))=f_{k-1}(\tr(F)).
\]
This implies that
\begin{equation}\label{TRACE}
\tr(\AA''_k)\subseteq \AA''_{k-1}.
\end{equation}

\section{Proof of Theorem \ref{AUTOMORPH}}

\subsection{A lemma about twin-free graphs}\label{TWINLEMMA}

We start with a simple lemma about twin-free weighted graphs.

\begin{lemma}\label{TWINFREE}
Let $G$ be a twin-free weighted graph. Then every map $\phi:~V(G)\to
V(G)$ such that $\beta_{\phi(i)\phi(j)}=\beta_{ij}$ for every $i,j\in
V(G)$ is bijective.
\end{lemma}

\begin{proof}
The mapping $\beta$ has some power $\gamma=\beta^s$ that is
idempotent. We claim that $i$ and $\gamma(i)$ are twins. Indeed,
\[
\beta_{il}=\beta_{\gamma(i)\gamma(l)}
=\beta{\gamma^2(i)\gamma(l)}=\beta_{\gamma(i)l}
\]
for every $l\in V(G)$. Since $G$ is twin-free, this implies that
$\gamma$ is the identity, and so $\beta$ must be bijective.
\end{proof}

\subsection{From Lemma \ref{HOMEQ} to Lemma \ref{HOMREP}
to Theorem \ref{AUTOMORPH}}\label{REDUCTIONS}

Now we turn to the proof of Theorem \ref{AUTOMORPH}. Let $\AA'_k$ be
the subalgebra of elements of $\AA_k$ invariant under the
automorphisms of $G$, and let $\AA''_k=f_k(\GG_k)$. It is trivial
that $\AA''_k\subseteq \AA'_k$. Furthermore, we have
\[
\dim(\AA'_k)=\frac{|V(G)|^k}{\orb_k(G)}
\]
and
\[
\dim(\AA''_k)=\dim(\GG'_k)=\rk(M(k,G)).
\]
Thus it follows that $\rk(M(k,G))\le |V(G)|^k/\orb_k(G)$; to prove
Theorem \ref{AUTOMORPH}, it suffices to prove that algebras $\AA'_k$
and $\AA'_k$ are the same. This is just the content of lemma
\ref{HOMREP}. Thus it suffices to prove this lemma.

The algebra $\AA''_k$ is a finite dimensional commutative algebra
with a unit element, and so it has a basis $w_1,\dots,w_r$ consisting
of idempotents. Expressing these idempotents in the basis of the
whole algebra $\AA_k$, we get that for each $i$ there is a set
$\Psi_i\subseteq V(G)^k$ such that
\[
w_i=\sum_{\psi\in\Psi_i} \psi.
\]
Since $\sum_k w_k$ is the unit element, it follows that the sets
$\Psi_i$ ($i=1,\dots,r$) form a partition of $V(G)^k$. We say that
$\phi,\psi\in[1,m]^k$ are {\it equivalent}, if they belong to the
same set $\Psi_i$. Clearly $\phi$ and $\psi$ are equivalent if an
only if
\[
\hom_{\phi}(F,G)= \hom_{\psi}(F,G)
\]
for every $k$-labeled graph $F$. The subalgebra $\AA''_k$ consists of
those elements in which any two maps $\phi$ and $\psi$ that are
equivalent occur with the same coefficient. Analogously, the
subalgebra $\AA''_k$ consists of those elements in which any two maps
$\phi$ and $\psi$ such that $\psi=\phi\sigma$ for some automorphism
$\sigma$ occur with the same coefficient. The fact that these two are
the same is just the content of lemma \ref{HOMEQ}. So it suffices to
prove this Lemma.

\subsection{Proof of Lemma \ref{HOMEQ}}\label{HOMEQPROOF}

For any map $\phi:[1,k]\to[1,m]$, let $\phi'$ denote its restriction
to $[1,k-1]$.

\begin{claim}\label{TRUNC}
If the maps $\phi,\psi\in[1,m]^k$ are equivalent, then so are $\phi'$
and $\psi'$.
\end{claim}

Indeed, assume that $\phi'$ and $\psi'$ are not equivalent, then
there is a $(k-1)$-labeled graph $F$ such that
\[
\hom_{\phi'}(F,G)\not= \hom_{\psi'}(F,G).
\]
Then for $F'=F\otimes E_1$ we have
\[
\hom_{\phi}(F',G) = \hom_{\phi'}(F,G)\not= \hom_{\psi'}(F,G) =
\hom_{\psi}(F',G),
\]
which contradicts the assumption that $\phi$ and $\psi$ are
equivalent.

\begin{claim}\label{EXTEND}
Suppose that $\phi,\psi\in[1,m]^k$ are equivalent. Then for every
$\mu\in[1,m]^{k+1}$ such that $\phi=\mu'$ there exists a
$\nu\in[1,m]^{k+1}$ such that $\psi=\nu'$ and $\mu$ and $\nu$ are
equivalent.
\end{claim}

Let $\Psi$ be the set of maps equivalent to $\mu$. By definition, we
have
\[
\sum_{\eta\in\Psi} \eta \in\AA''_{k+1}.
\]
Applying the trace operator, we see by (\ref{TRACE}) that
\[
\sum_{\eta\in\Psi} \alpha(\eta(k+1)) \eta' \in\AA''_k.
\]
Here $\phi$ occurs with non-zero coefficient; since $\phi$ and $\psi$
are equivalent, $\psi$ must occur with non-zero coefficient, which
shows that there must be a $\nu\in\Psi$ such that $\nu'=\psi$. This
proves Claim \ref{EXTEND}.

To prove the Lemma, we want to show that if $\phi$ and $\psi$ are
equivalent, then there exists an automorphism $\sigma$ of $G$ such
that $\psi=\phi\sigma$. We prove this assertion for an increasing
class of mappings.

\begin{claim}\label{ONE2ONE}
If $k=m$ the maps $\phi,\psi:~[1,m]\to [1,m]^m$ are equivalent, and
$\phi$ is bijective, then there is an automorphism $\sigma$ of $G$
such that $\psi=\phi\sigma$.
\end{claim}

We may assume that the nodes of $G$ are labeled so that $\phi$ is the
identity. Let $\Psi$ be the set of maps equivalent to $\phi$. We
claim that every $\psi\in\Psi$, viewed as a map of $V(G)$ into
itself, satisfies
\begin{equation}\label{BETAPRES}
\beta_{ij}=\beta_{\psi(i)\psi(j)}
\end{equation}
for every $j$. Indeed, let $k_{ij}$ be the $k$-labeled graph
consisting of $k$ nodes and a single edge connecting nodes $i$ and
$j$. Then
\[
\beta_{ij}=\hom_{\phi}(k_{ij},G)
=\hom_{\psi}(k_{ij},G)=\beta_{\psi(i)\psi(j)}.
\]
Since $G$ is twin-free, it follows by Lemma \ref{TWINFREE} that
$\psi$ is one-to-one.

To complete the proof of the Claim, it suffices to show that for
every $\psi\in\Psi$,
\begin{equation}\label{APRES}
\alpha(j)=\alpha(\psi(j))\qquad (j=1,\dots,m).
\end{equation}
It suffices to prove this for the case $j=m$. By the definition of
equivalence, we have
\[
\sum_{\psi\in\Psi} \psi \in\AA''_m.
\]
Applying the trace operator, we see by (\ref{TRACE}) that
\[
\sum_{\psi\in\Psi} \alpha(\psi(m)) \psi' \in\AA''_{m-1}.
\]
As we have seen, all maps $\psi\in\Psi$ are bijective, which implies
that the maps $\psi'$ are all different. Since these maps are
equivalent by Claim \ref{TRUNC}, it follows that all coefficients
$\alpha(\psi(m))$ are the same. This completes the proof of Claim
\ref{ONE2ONE}.

\begin{claim}\label{ONTO}
If the maps $\phi,\psi:~[1,k]\to [1,m]$ are equivalent, and $\phi$ is
surjective, then there is an automorphism $\sigma$ of $G$ such that
$\psi=\phi\sigma$.
\end{claim}

By permuting the labels $1,\dots,k$ if necessary, we may assume that
$\phi(1)=1,\dots,\phi(m)=m$. Claim \ref{TRUNC} implies that the
restriction of $\psi$ to $[1,m]$ is equivalent to the restriction of
$\phi$ to $[1,m]$, and so by Claim \ref{ONE2ONE}, there is an
automorphism $\sigma$ of $G$ such that $\psi(i)=\sigma(i)$ for
$i=1,\dots,m$.

Consider any $m+1\le j\le k$, and let $\phi(j)=r$. We claim that
$\psi(j)=\psi(r)$. Indeed, the restriction of $\phi$ to
$\{1,\dots,r-1,r+1,\dots,m,j\}$ is bijective, and equivalent to the
restriction of $\psi$ to this set; hence the restriction of $\psi$ to
this set must be bijective, which implies that $\psi(j)=\psi(r)$.
This implies that for every $1\le i\le k$, $\psi(i)=\sigma(\phi(i))$
as claimed.

Now we are ready to prove the Lemma for arbitrary equivalent maps
$\phi,\psi:~[1,k]\to[1,m]$. We can extend $\phi$ to a mapping
$\mu:~[1,\ell]\to[1,m]$ ($\ell\ge k$) which is surjective. By Claim
\ref{EXTEND}, there is a mapping $\nu:~[1,\ell]\to[1,m]$ extending
$\psi$ such that $\mu$ and $\nu$ are equivalent. Then by Claim
\ref{ONTO}, there is an automorphism $\sigma$ of $G$ such that
$\nu=\mu\sigma$. Restricting this relation to $[1,k]$, the assertion
follows.

\subsection{Proof of Corollary \ref{HOMDET}}

Let $G$ be the graph obtained by taking the disjoint union of $G_1$
and $G_2$, creating two new nodes $v_1$ and $v_2$, and connecting
$v_i$ to all nodes of $G_i$. Also add loops at $v_i$. The new nodes
and new edges have weight 1.

We claim that for every 1-labeled graph $F$,
\begin{equation}\label{V1V2}
\hom_{v_1}(F,G)=\hom_{v_2}(F,G)
\end{equation}
Indeed, if $F$ is not connected, then those components not containing
the labeled node contribute the same factors to both sides. So it
suffices to consider the case when $F$ is connected. Then we have
\[
\hom_{v_1}(F,G)=\sum_{S\subseteq V(F)\atop S\ni v_1} \hom(F\setminus
S,G_1).
\]
Indeed, every map $\phi:~[1,k]\to V(G)$ such that $\phi(1)=v_1$ maps
some subset $S\subseteq V(F)$, $S\ni v_1$ to the new node $v_1$; if
we fix this set, then the restriction $\phi'$ of $\phi$ to
$V(F)\setminus S$ is a map into $V(G_1)$ (else, the contribution of
the map to $\hom(F,G)$ is $0$), and the contribution of $\phi$ to
$\hom_{v_1}(F,G)$ is the same as the contribution of $\phi'$ to
$\hom(F,G)$.

Since $\hom_{v_2}(F,G)$ can be expressed by a similar formula, and
the sums on the right hand sides are equal by hypothesis, this proves
(\ref{V1V2}).

Now Lemma \ref{HOMEQ} implies that there is an automorphism of $G$
mapping $v_1$ to $v_2$. This automorphism gives an isomorphism
between $G_1$ and $G_2$.

\section*{Acknowledgement}

I am indebted to Christian Borgs, Jennifer Chayes, Mike Freedman,
Monique Laurent, Lex Schrijver, Miki Simonovits, Vera T.~S\'os,
Bal\'azs Szegedy, G\'abor Tardos and Kati Vesztergombi for many
valuable discussions and suggestions on the topic of graph
homomorphisms.


\begin{thebibliography}{99}

\bibitem{FLS}
M.~Freedman, L.~Lov\'asz and A.~Schrijver: Reflection positivity,
rank connectivity, and homomorphism of graphs (manuscript)

\bibitem{Lo}
L.~Lov\'asz: Operations with structures, {\it Acta Math.\ Hung.} {\bf
18}, 321-328.

\bibitem{LS}
L.~Lov\'asz and V.T.~S\'os: Generalized quasirandom graphs,
manuscript.
\end{thebibliography}
\end{document}